\documentclass[12pt]{amsart}

\usepackage{amsmath}
\usepackage{amssymb}
\usepackage{url}
\usepackage[greek,english]{babel}

\usepackage{color}
\def\sideremark#1{\ifvmode\leavevmode\fi\vadjust{
\vbox to0pt{\hbox to 0pt{\hskip\hsize\hskip1em
\vbox{\hsize1cm\tiny\raggedright\pretolerance10000
\noindent #1\hfill}\hss}\vbox to8pt{\vfil}\vss}}}

\def\hl#1{#1}

\large
\title{Geometry in the Transition from Primary to Post-Primary}
\author{Patrick D. Barry and Anthony G. O'Farrell}
\address[P.D. Barry]{School of Mathematical Sciences,
UCC, Cork, Ireland}
\address[A.G. O'Farrell]{Department of Mathematics and Statistics, 
NUI,\break\null{}\qquad
Maynooth, Co. Kildare, Ireland}
\email[Barry]{p.barry@ucc.ie}
\email[O'Farrell]{admin@maths.nuim.ie}

\begin{document}
\maketitle
\section{Introduction}

This article is intended as a kind of precursor to the 
document {\it Geometry for Post-primary
School Mathematics}, which forms Section B, pp. 37--84
of the Mathematics Syllabus for Junior Certicate
issued by the National Council for Curriculum
and Assessment \cite{NCCA} in the context of 
Project Maths.

Our purpose is to place that document in the context 
of an overview of plane geometry, touching on several 
important pedagogical and historical aspects, in the
hope that this will prove useful for teachers.

The main points we want to emphasize are these:

\begin{itemize}
\item Geometry is a key part of mathematics.
\item Children must pass through different stages in studying 
geometry.
\item Each stage plays an important r\^ole.
\item Care must be taken in managing the transitions.
\item Some knowledge of history is useful for teachers.
\end{itemize}

\section{The Main Parts of Mathematics}
At present, the NCCA presents the mathematics 
curriculum in terms of {\em strands}.
For primary level the five strands are labelled
Number, Algebra, Shape and Space, Measures,
and Data. For secondary level they are 
(1) data, statistics and probability, 
(2) geometry and trigonometry,
(3) number and measure, (4) algebra and (5) functions.
The similarity between the two
classifications  is part of an attempt to
foster continuity across the transition. A further
initiative designed to foster this was the publication
of a {\em bridging framework} \cite{bridging}
which provides a dictionary linking terminology used in
primary to that used in secondary schools.
Initially, the divisions were perhaps intended as
much to reflect a more-or-less equal
division of teaching and learning effort as much as
a division of mathematics into its main areas.
It was never intended that the strands would represent
watertight divisions of the subject, and was generally
recognised that there is necessary interaction between 
them.

An analysis of the secondary curriculum reveals that
the division becomes progressively more forced
at the more advanced levels.  It would also be
completely impossible to impose it on third-level
studies in mathematics.  It is standard among
mainstream professional researchers to say that
there are three main branches in the tree of
mathematics: Algebra, Analysis and Geometry.
Some would add Probability and Computation
to these, but many others would regard 
Probability and Computation as two branches of 
Analysis.
However, is is also standard that the connections between
the main fields of mathematics are so many that
it is actually possible to say that any one of these
embraces the whole.  There are fields such as 
Algebraic Geometry, Algebraic Topology and
Geometric Analysis which may be regarded as
branches of either one of two main branches, and which
use fundamental results from both. 

The point we make here is that from a bird's-eye
viewpoint, Geometry is about one third of mathematics.
At the research level, it accounts for a solid
proportion of new PhD theses, as may be seen by 
examining the tables in reports of the AMS annual surveys 
of new graduates
\cite{survey}.  
Its applications
are in active areas of fundamental and applied physics,
robotics, coding, graphics and other commercially 
significant areas.  So it is important.  Our
students deserve a sound formation in geometry.

For various reasons, in the past decades many students have 
emerged from secondary school with a poor opinion of geometry,
the result of unfortunate experience with the subject.  
PISA assessment results also showed relatively mediocre performance by
the Irish 15-year-old cohort on problems requiring
geometrical skill.
This needs to change, and indeed change is mandated by
the Project Maths curriculum.

\section{Stages}
\subsection{Primary stages}
The present Primary Curriculum \cite{primary}
specifies the study of geometrical shapes
in two and three dimensions under 
the heading of Shape and Space, and of length
and area under Measure (with support from Number).
This begins right at the start, and
is developed further year by year.
Students are introduced to simple
shapes (triangle, rectangle, circle, semicircle, cube,
cuboid, sphere, cylinder, cone) and progressively more
complex shapes and properties (isosceles triangle,
parallelogram, rhombus, pentagon, hexagon,
triangular prism, pyramid, scalene triangle,
trapezium, regular hexagon)
learn to distinguish them and learn names for them
and for their parts and properties.
They make use of suitable materials (blocks, paper and scissors,
folded paper,
art straws, geoboard, mazes, grids, board games,
software, plasticine, prisms, compass, string,
tangrams, squared paper) and diagrams
and learn to recognize shapes in their environment.
Most of the work involves flat, planar shapes,
but they also manipulate 3-D shapes and solve problems
about them.
They learn about measuring lengths, areas, volumes
and angles (using a progressively
richer number system).  They learn how the describe and
evaluate spatial relations, give directions,
construct and draw 2-D shapes using instruments, 
subdivide and combine shapes, draw tesselations,
construct 3-D polyhedra (by folding nets), 
and use coordinates.
They are encouraged to look
for common patterns such as lines of symmetry
and the result of counting faces minus
edges plus vertices
for polyhedral shapes.  
They learn about parallel lines, and
right, acute, obtuse
and reflex angles.
They explore properties of 2-D shapes,
including the angle-sum of a triangle
and a quadrilateral, and the ratio
of diameter to radius of a circle. 
The syllabus specifies the linkage of
geometry to other areas of the curriculum
(motor skills, science, art, physical education and dance,
geography) and to aspects of everyday life.
This is all very useful, and is appropriate for their
ages.  

At the end of primary school, children should have
acquired most basic geometrical concepts and the language
that goes with them.  They should be in a position
to use their understanding to solve many practical
problems.

\medskip 
\subsection{Secondary stages}
When they start post-primary school, students should not be allowed to abandon 
all this geometrical experience, but should continue to draw on it,
solidify and develop their understanding of it, and 
stay in touch with geometrical ideas on a continual basis. 

They have a lot more to learn about geometry.
There are in fact two further components needed beyond
primary level, corresponding to the
two main reasons that further geometrical study is needed:
the practical utility of more advanced material and skills,
and the r\^ole of geometry in developing and honing
the student's reasoning power.

The case for exposure to rigorous mathematical
thought as a preparation for life and for any
further studies was well made by
John Stuart Mill (quoted in \cite{Potts}):

\begin{quotation}

{\it{ }
The value of Mathematical instruction as a preparation
for those more difficult investigations (physiology,
society, government, \&c.) consists in the applicability
not of its doctrines, but of its method.  Mathematics 
will ever remain the most perfect type of the Deductive
Method in general $\cdots$ 

These grounds are quite sufficient for deeming
mathematical training an indispensible basis of real
scientific education, and regarding, with Plato,
one who is\\ \greektext `agewm'etrhtoc%
\footnote{\latintext -- {\it ageometretos }, i.e. ignorant of geometry
(or, perhaps, unskilled in geometry, or indifferent to geometry).  
The motto said to have been
carved above the entrance to Plato's
Academy was: \greektext O`udeic `agewm'etrhtoc e`is'itw --
\latintext {\it Let no-one
ignorant of geometry enter.}}
\latintext as wanting in one of the
most essential qualifications for the successful cultivation
of the highest branches of philosophy.
}
\end{quotation}  

Geometry is not the only branch of Mathematics
that may serve to develop reasoning power,
but it has long been used for that purpose,
and many consider it well-suited.  The geometrical
theory expounded in the Elements
of Euclid (cf. \cite{Euclid}), rediscovered in the
West at the end of the Middle Ages and
adopted as the preferred text
by the first European universities has
been the most popular.  It is important for
teachers to understand some key points about it:

\begin{itemize}
\item It is an abstract theory about space (without matter).
\item It was not written to be studied by children.
\item It has some logical flaws.
\item Mathematicians have figured out various ways to fix these
flaws so that the main propositions can be proven
from a set of axioms.  
\item Each such amended theory is called Euclidean Geometry.
\item Each such theory is even less suited for children.
\item Euclidean geometry is very useful.
\item There are other geometrical theories, in which
some of the propostions of Euclidean Geometry are false.
\item We do not actually know which of these is the
best approximation to \lq\lq real empty space"\footnote{
-- another abstraction.}.
\item We do know that actual space, with matter,
does not fit well with Euclidean Geometry\footnote{when lines are interpreted
as light rays.},
although there is a close correspondence at
everyday length scales.
\item Abstract geometry has to be simplified, if it is to be used 
in school to develop reasoning power.
\item Even when simplified, it is not feasible just to
fling children into the abstractions without a careful
preparatory stage.
\end{itemize}

We shall elaborate on some of these points,
and comment on the pedagogical implications.

We start with some history.

\section{Historical development of geometrical theory}
\subsection{The arc of history}
Euclidean Geometry has had a long history. Following on practical studies 
of shapes, lengths, areas and volumes in the Sumerian and Egyptian 
civilizations \textit{inter alia,}  it started to evolve into a 
logically-organised science as a result of
the efforts of philosophers in Greece c.700-600 BC, 
who wanted to base knowledge on solid foundations. 
The basic idea was to identify and define
purely geometrical (i.e. non-material) abstractions
(point, line, etc.)
and also identify uncontroversial starting principles
about them, and then to use logic to work out
the consequences.
Ideally, the building blocks of the theory should be
as simple as possible.
This simple idea proved extremely effective in
practical applications (such as tunnel construction), 
and gave encouragement.
Understanding of geometrical theory
evolved gradually ever since, although there
were  many fallow centuries. 
Euclid's synthesis of the geometry of his day (about 300 BC)
was a
major landmark, but after his time many further
theorems unknown to him have been
discovered, and our understanding of the
basic plan of his work has also evolved.
Ren\'e Descartes, in his Discourse on the Method (1637 AD) showed
how to link numbers to geometry -- in Euclid's books,
and up to that point
geometrical magnitudes and numerical magnitudes
had been considered different species.  This created
the field of algebraic geometry, and in a sense reduced
geometry to arithmetic.  However 
it may come as a surprise that the system of real numbers has been 
fully understood only since about 1860 AD (thanks to Richard Dedekind).
Non-euclidean geometries were discovered early in the 19th
century, rubbishing Kant's view that our knowledge of
Euclidean geometry is \lq\lq synthetic a priori",
and raising the question whether the real world is
Euclidean or not.  

\subsection{Deductive reasoning}
A proper understanding of logical deductive 
systems was only arrived at in the late 1800's,
and this prompted Hilbert to produce the first
fully-rigorous account of Euclidean geometry,
i.e. an account in which all the theorems of
Euclid can be proved rigorously from first principles. 
What is now understood as a mathematical theory,
or deductive system, has five components \cite{Greenberg}:

\begin{enumerate}
\item Undefined terms 
\item Definitions
\item Axioms
\item A system of logic (rules for valid deductions)
\item Theorems (a term that embraces also propositions, lemmas, and
corollaries).
\end{enumerate}

In other words,
in a logical system we list up-front the terms
and assumptions that we start with, and thereafter proceed by way of
definitions and proofs. 

\subsection{Definitions}
Definitions are about specifying what we are dealing with. A
\textit{definition} identifies a new concept in terms of accepted or known
concepts. In practice a definition of a word, symbol or phrase $E$ is a
statement that $E$ is to be used as a substitute for $F$, the latter
being a phrase consisting of words and possibly symbols or a compound
symbol. We accept ordinary words of the English language in definitions
and what is at issue is the meaning of technical mathematical words or
phrases. In attempting a definition, there is no progress if the
technical words or symbols in $F$ are not all understood at the time of
the definition.

The disconcerting feature of this situation is that in any one
presentation of a topic there must be a first definition and of its
nature that must be in terms of accepted concepts. Thus we must have
terms which are accepted without definition, that is there must be
\textit{undefined} or \textit{primitive} terms. This might seem to leave us in a
hopeless position but it does not, as we are able to, and must, assume properties
of the primitive terms and work with those.

There is nothing absolute about this process, as a term which is taken as
primitive in one presentation of a topic can very well be a defined term
in another presentation of that topic, and vice versa. We need
\textit{some} primitive terms to get an approach under way.

\subsection{Proof}
Proof is the way to establish
the properties of the concepts that we are
dealing with. A \textit{proof} is a finite sequence of statements the first
of which is called the \textit{hypothesis}, and the last of which is called
the \textit{conclusion}. In this sequence, each statement after the
hypothesis must follow logically from one or more statements that have
been previously accepted. Logically there would be a vicious circle if
the conclusion were used to help establish any statement in the proof.

There is also a disconcerting feature of this, as in any presentation of
a topic there must be a first proof. That first proof must be based on
some statements which are not proved (at least the hypothesis), which
are in fact properties that are accepted without proof. Thus any
presentation of a topic must contain unproved statements; these are
called \textit{axioms} or \textit{postulates} and these names are used
interchangeably.

Again there is nothing absolute about this, as properties which are
taken as axiomatic in one presentation of a topic may be proved in
another presentation, and vice versa. But we must have {\em some} axioms
to get an approach under way.

\subsection{Hilbert's system}
In Hilbert's system \cite{Hilbert} there
are undefined terms such as {\em point, line,
plane, between, congruent}, seven axioms of
connection, five axioms of order, an axiom of
parallels, six axioms of congruence, and
an axiom of continuity, and definitions of terms such as 
{\em segment, vertex, side of a line}.
The logic used is standard Aristotelian logic.

Notice that this leaves aside completely the question
of any relation between this theory and the real world.
There are equally satisfactory 
and equally-consistent\footnote{The consistency 
of Euclidean geometry cannot be proven.  It can be
shown that it is consistent if elementary
arithmetic is consistent.} theories of various geometries
in which some of Euclid's theorems are false.

Incidentally, the main aspects of Euclid's work
that needed to be \lq\lq cleaned up" were
(1) the attempt to prove the SAS
congruence criterion, Prop. I:4,
using superposition,
instead of just assuming it;
(2) the absence of any postulates about
line separation or plane separation,
and \lq\lq betweenness", needed for instance
in Prop. I:16; and
(3) the absence of any continuity or
completeness assumption, already an
issue in Prop I:1.
The notion gained currency in the 1960's
that \lq\lq Euclid is all wrong" and
should just be dumped.  The truth is that
with a little careful tweaking early on,
almost everything stands and the proofs can
still be used.

\subsection{Alternative Versions of Euclidean Geometry}

Over the period from c.500 BC to the present 
quite a few different 
approaches to Euclid's theorems
have been published. The superabundance of these is 
one of the major problems that we face now. 
Hilbert's was the first that was logically watertight and
categorical\footnote{-- in the sense that any two
interpretations (models) of it are essentially the same}.
Before his time, and since, many variants were invented by
teachers who wanted to make Euclid more accessible
to children.  After him, other professional 
research mathematicians produced complete versions
involving different undefined terms, definitions and
axioms, but of course the same theorems.  They were 
motivated by the desire to have an equivalent system
with simpler axioms.  For instance, 
Hilbert's system does not
include or use the real number system, and Birkhoff
\cite{Birkhoff} proposed a system that extended the theory of the
real numbers by adding only four axioms and gave
all of Euclidean geometry.  

\section{Geometry for secondary school}
\hl{
\subsection{}
It should be 
stated honestly, and faced now that a fully rigorous 
account of plane 
Euclidean geometry is too difficult for secondary school.
This has been widely understood in academic circles
for a very long time.  In Mathematics Education circles,
it was made explicit in the work of the van Hieles (cf. \cite{Crowley})
in the nineteen-fifties, when they identified five levels
at which a person might understand geometry, ranging from level
zero (\lq\lq visualization") up to level four (\lq\lq rigour").
The top level is only appropriate for university-level work.
This does not mean that logical work in geometry is not
feasible in school --
the van Hiele levels appropriate for school were labelled
\lq\lq analysis", \lq\lq informal deduction" and
\lq\lq formal deduction".  Moreover, competence at
the top level is not 
really needed for working with the manifold applications
of geometry\footnote{%
The van Hiele model has had a lot of influence.
Early on, it formed the basis for a radical reform in the
geometry curriculum in the Soviet Union in the 
nineteen-sixties, and it has gradually been taken
on board in the USA.}. 
}

\hl{
A key point is that you cannot train someone in logical
deductive thinking by using an illogical system.
So what professional 
mathematicians urge and press for is that school geometry should be in the 
the context of some fully-rigorous scholarly background 
approach. The school version should broadly 
have the same sequence of topics and the same type of proofs, but 
leave out some very difficult parts, the latter to be guaranteed 
by those at a higher level who choose to immerse themselves in a study 
of this material.
}

The present situation in Ireland is that the prescribed
school geometry has for its scholarly background approach
the one laid out in \cite{Barry}.  That system, like Birkhoff's,
uses the real number system, but employs a few more axioms,
including Playfair's version
of the parallel axiom \cite[Axiom A7, p.57]{Barry}.
The school system is deliberately simplified, as explained in
\cite[pp.40-43]{NCCA}.  As a result, the proofs are not
fully watertight, relying in places on unstated
\lq\lq commonsense notions", and a teacher or student
who notices this is encouraged to refer to \cite{Barry}
to satisfy themselves that the gaps can be bridged.

An important aspect of Birkhoff-like systems is that
one can treat the real number system informally
(instead of formally, axiomatically)
in school.  This avoids explicitly mentioning
the topics of continuity and completeness, which
are too sophisticated for school.  An important reason for
basing the system on the book \cite{Barry}
is that there the complete scholarly treatment is
fully laid out, {\em with complete detailed proofs}. A number
of good alternatives are backed by complete theories
for which the full proofs can be generated easily by
any competent professional mathematical researcher,
but are in fact only sketched in published sources.

\subsection{Approaching abstraction when teaching geometry}

Euclidean geometry employs abstractions.
Right at the start, we have point and line, for instance.

Students have to be prepared carefully for this abstraction.
The geometrical concepts must be motivated 
from the real world around us. 
Education in geometry (as in everything else) must
proceed in stages, as the child's mind develops.
These stages have long been recognised, and were explicitly
catered for in popular textbooks such as Durrell
\cite{Durell}.  (Clement Durell's texts were in widespread
use wherever British influence acted
from 1919 for over forty years.)
There has to be a preliminary stage before the
stage of formal logical work with the abstract ideas
and it is essential that these stages not be confused with each other. 
The preliminary stage should not be rushed,
and time allowed for the abstract concept to sink in.
It is not appropriate to plunge into \lq\lq Theorem 1"
immediately after explaining about points. 

Later, when abstract
results are applied, we should make it clear that
we are now assuming they apply to reality.

A point is not a real thing.  It has no size.  
Durell \cite{Durell} says
that teachers should {\em never} allow points to be 
drawn as blobs, and instead indicated by a cross made with
two very fine lines. He insists that compass punctures
should be as small as possible, and straight lines be as
fine as possible. This is extreme, but
you see his point!  There may well be students who think
points are little black round things, as drawn by
Geogebra, and it is a good idea to make sure that
they are disabused of this before they get started
on formal work. 
 
Diagrams are vital in teaching geometry, and should always
be used. It is precisely
because such visual aids are there to support and guide reason
that geometry is considered the best way to
practise logic.

There should be considerable physical motivation to 
start with, and diagrams always used to provide insight, but 
details of motivation should not be confused with the careful 
logical presentation of the mathematical model that 
follows later on.

Every opportunity should be taken to get students to engage with
problems that they can tackle using their current understanding
of geometry.

\hl{
As with any mathematics teaching, one proceeds
in a cycle \cite{Durell, Crowley}: oral discussion of examples,
exercises in numerical and non-numerical examples, 
informal proof ideas, formal proofs,
}
exercises involving \lq\lq riders", or \lq\lq cuts"
(extra propositions to be proven by the student
-- the \lq\lq Propositions" given without proof
in the syllabus document are intended to be
used in this way, and it is expected that the
assessment process will examine skill in
creating such proofs),
and one provides exercises graduated by difficulty,
extra exercises of one kind for students who struggle, and
challenging extra exercises for those pupils \lq\lq who
run ahead of the class".
Regarding the latter, although the main focus of the school 
programme is on plane geometry, one should look out for
applications to solid geometry.

\section{Lines and Non-Euclidean Geometry}
The modern mathematical concept of {\em line} is infinite, without ends,
and is straight. The English word line is derived from
the Latin {\it linea}, which originally referred to
flaxen thread, as is the name of the material {\em linen},
also made from flax. Similarly, in Irish we have
the pair of words \textit{l\'{i}ne} and \textit{l\'{i}n\'{e}adach}. 
In Greek, the word for line was \greektext gramm'h \latintext (gramme), 
the stroke
of a pen, derived from \greektext grafw \latintext (grapho), I write, or draw.
In contrast to modern usage 
the Greeks spoke of a \textit{straight line} (literally
\greektext e`uje\~ia gramm'h \latintext, \lq\lq right line") 
and \textit{curved line}.
Moreover, by straight lines the Greeks mainly meant what we
call \textit{line-segments} which would be \textit{produced} (i.e. extended)
as required. It is helpful to bear this in mind when reading
older texts. 

One motivation of a line-segment was a \textit{linen 
thread held taut}\footnote{Heath\cite{Euclid} records the classical phrase
\greektext 'ep >akron tetam'enh gramm'h \latintext
-- a line stretched to the utmost.}. 
The notion of being straight was extended to lines, 
as a segment was unendingly produced, 
and at each stage there had to be a segment which contained the 
starting segment. Of course a taut thread could be copied onto a 
wax, papyrus or wooden tablet, and tablets with straight edges 
could be cut from the latter. The use of compasses then enabled 
them to cut out shapes of triangles and various types of quadrilaterals, 
as well of course as circles. Nowadays we have rulers, protractors, 
set squares and computer software packages to help us draw these 
figures and make constructions.

Euclid defined a straight line as \lq\lq a line that lies evenly
with the points on itself".  This suggests that the concept is an 
abstraction of the idea of a \lq\lq line of sight", or if you
prefer, of a light ray, although Euclid in the Elements
is careful not to refer to any physical thing, perhaps because
as a good Platonist he is treating geometry as the form
of space, a purer thing than space (cf. \cite[pp. 165-169]{Euclid}).

If straight lines are the paths of light rays from an object
to our eye, then we now know, from observational Astronomy,
that real space is not Euclidean: there are distant objects that
can be seen in two different directions ( -- usually 
explained in terms of the theory of General Relativity as
the result of \lq\lq gravitational lensing").  
So the most that might be true is that real space
is approximately Euclidean at some length-scales
and where the density of matter is low.  
Gauss realized that space might be non-Euclidean
before most people, and checked measurements
taken of the angle-sum in large terrestrial triangles
(with vertices on mountain-tops in the Harz) as part of his
geodetic survey of the kingdom of Hanover
in the early 1820's.  These measurements did not
show a deviation from $180^\circ$, within experimental
error.  In certain non-Euclidean geometries ("hyperbolic")
the angle-sum is always less than $180^\circ$, and the defect
grows larger and larger as the triangle gets larger --
in extremely large triangles the sum may be arbitrarily
small!  The defect is small in small triangles\footnote{
In other non-Euclidean geometries (\lq\lq elliptic", or \lq\lq spherical")
the angle-sum is
consistently greater than $180^\circ$}.
It may be that if we could measure the angles
in a triangle with vertices in three different galaxies
( -- hard to see how to do this without visiting them),
we would find it substantially less than $180^\circ$.
We don't know.
So we are not justified in instructing our students
to \lq\lq discover that the angle-sum is $180^\circ$".
They can and should discover that it appears to be very
close to this, but that is all. 
In a class of students, measuring hand-drawn triangles
to half a degree,
one might expect a range of measurements clustering around $180^\circ$,
and that is fine.  

The proof of the angle-sum theorem, Euclid I:32,
\cite[Theorem 4]{NCCA}, just shows that if you accept the Axioms,
then the theorem holds.  It is probably not advisable
to disturb the faith of the young, but it is as well
for teachers to understand this.  It hinges on the Axiom
of Parallels, and it is amusing to list some of the
alternative assumptions that could be used instead
(along with the other axioms\footnote{-- strictly, the
other axioms of \lq\lq neutral geometry" \cite{Greenberg}}) , and
that are false if it fails (the names in brackets
after each are associated with them):

\begin{itemize}
\item There are two similar (equiangular) triangles
in which each side of one is twice the
corresponding side of the other (Wallis, 1663)
\item There is at least one rectangle (Saccheri, 1773; Omar Khayyam, 11th Century).
\item There is a triangle having area as
great as you please (Gauss, 1799).
\item A line perpendicular to the bisector of an acute angle
at a point inside the angle must meet both arms of the
angle. 
\item A line that cuts one of two parallel lines must cut the other 
(Proclus, 5th Century AD).
\end{itemize}

An internet search for the parallel postulate will throw up
many more such oddities.  The fact that there are
so many very plausible statements that imply the Parallel
Postulate explains why so many eminent mathematicians
were deceived into thinking they had proved it without
assuming anything\footnote{Courses in non-Euclidean geometry
are quite usual in the preparation of mathematics teachers
these days. To quote Wolfe \cite{Wolfe}:
\lq\lq For teachers and prospective teachers of geometry in the secondary 
schools the study of Non-Euclidean Geometry is invaluable. 
Without it there is strong likelihood that they will not understand the 
real nature of the subject they are teaching and the import of its 
applications to the interpretation of physical space."
He in turn quotes Chrystal, who published a small book
about what he called pan-geometry in 1880, aimed at teachers.
He wrote:
\lq\lq It will not be supposed that I advocate the introduction of 
pan-geometry as a school subject; it is for the teacher that I advocate 
such a study. It is a great mistake to suppose that it is sufficient for 
the teacher of an elementary subject to be just ahead of his pupils. 
No one can be a good elementary teacher who cannot handle his 
subject with the grasp of a master. Geometrical insight and wealth 
of geometrical ideas, either natural or acquired, are essential to a 
good teacher of geometry; and I know of no better way of 
cultivating them than by studying pan-geometry."  
}. 

This explains why there is no such thing
as a \lq\lq physical proof" of a mathematical
theorem, and why no theorem is \lq\lq visually
obvious".  Any attempt to deduce a mathematical
fact from a real-world observation involves
a logical gap: the implicit assumption that some
mathematical theory accurately reflects some
physical reality.  If such an assumption is
made explicit, it becomes a scientific hypothesis,
and can only be disproved by observation, never
confirmed \cite{Popper}. 

\medskip 
It does not seem to be widely appreciated that,
from the logical point of view, the abstract results
are also needed in order to lay a firm ground
for trigonometry and for coordinate geometry.
If the parallel postulate fails, we have no
rectangles, hence no rectangular cartesian coordinates, and
triangles are not similar unless congruent,
so we have no standard trigonometry.
The theorems about ratios have recently been
restored to the Leaving Certificate program,
in part because they provide this foundation
(retrospectively) for Junior Certificate 
trigonometry and coordinate geometry, and in part
because it was considered important that senior
students continue to engage with formal proofs
in synthetic geometry.

There may be nothing in nature that corresponds exactly to  
Euclidean geometry, 
but it cannot be denied that it has been 
extraordinarily useful in practical matters for over two and a 
half millennia.
It should also be noted that even if it does not
fit exactly the shape of the real universe,
Euclidean space, {\em as an ideal mental construct invented
by us}, is immensely useful in other areas of pure and
applied mathematics\footnote{Examples are the theory
of equations, numerical computation, much of 
real and complex analysis, and even non-Euclidean
geometry, which is studied using coordinates (\lq\lq charts")
that map pieces of the space to Euclidean space.}
and will always be used. By the way, the useful software system Geogebra
is a realization of this ideal mental construct. In it
the angle sum is always {\em exactly} $180^\circ$!

\section{More History: a fork in the road}

We now say a bit about the history of education in geometry.

\subsection{} Euclid's \textit{Elements} had a virtual monopoly as a textbook for geometry for a very long time. A substantial splinter-group was started in France in the 16th century when Pierre de la Ram\'{e}e (in Latin Petrus Ramus)(1515-1572), among his other publications, attacked the logical approach of the \textit{Elements}. His views attracted considerable support in French educational circles for many generations and led to widely-held views that what is visually obvious should be accepted without proof.  All this led to quite a different approach to geometry with many innovations. Some of these were reflex angles and rotations. Many of these ideas accumulated eventually to efficient new approaches to geometry, in modern times e.g to treatments based on transformations on the one hand and to vector spaces on the other. Thus there are now several possible approaches to Euclidean geometry available but, for present purposes, consideration of them should not be confined to their abstract merits but to which are the most suitable for our school students to obtain their grounding in geometry.

\subsection{} 
The lead of France was generally followed on the continent of Europe and the notable large country standing out against this and adhering to Euclid's approach was Great Britain. There there were efforts from 1860 onwards to assimilate elements then current on the Continent. The final bulwark to modifying Euclid's \textit{Elements} fell in 1903
with recommendations of the Mathematical Association,
and the \lq Cambridge Schedule'
proposed and adopted at the University of Cambridge about that time.
Subsequently in textbooks there, and used in Ireland, a variety of approaches and concepts were mingled from different strands \cite{L,MA,HS,NG}

\subsection{} 
This intermingling leads to some severe technical problems. All the textbooks 
started roughly the same way, focusing on the concepts in the world about 
us, becoming familiar with shapes and sizes and being led to properties of 
classes of them. On the whole they were clinging to a bad old habit from 
Euclid of trying to define everything.
  
From what we have said above it should be clear that 
they should have been motivating concepts from the real world, but that in 
formal geometry undefined but named items are needed. From Euclid they 
retained the concepts of assumed axioms, or postulates.

The practical difficulty in this is that for some concepts which are 
obvious and readily understood visually, it is quite difficult to lay 
down assumed properties for them which lead to their being singled out 
uniquely. In the next section we deal with two concepts which are at the heart of this problem.

By the way, it is a mistake to think that our forebears
were not aware of the need for undefined terms.  
In books written for
undergraduates, as opposed to school children,
they express themselves frankly on the point. When
Maynooth College was set up in 1795, it was initially staffed
by French clerical academics, refugees from the
Revolution. Mathematical instruction was compulsory for
all students, lay and clerical,
and it appears that geometrical teaching
was based on French models
so that  practice at Maynooth was in step with
the continental, rather than the British norm.
Andr\'e Darr\'e, first Professor of Mathematics
and Natural Philosophy, and formerly of Toulouse,
prepared a text in English \cite{Darre} on plane
and spherical geometry for use in the College.
He first gives the following definitions:

\begin{quotation}
{\it A straight line is that of which the elementary parts run
in the same direction. A line is curved, the elements
of which change continually their direction.}
\end{quotation}

Then he says:

\begin{quotation}
{\it Such is the most accurate notion Geometry can give of its
object; and it is adequate to its object, though not
perhaps a logical definition.

Sciences mostly begin by such simple ideas sufficiently clear,
independently of a definition; and they are no less reasonable than
self-evident.

For want of such simple notions and self-evident principles,
an interminable series of definitions and demonstrations 
should be required; our mind could find no\\ ground
whereon to rest in analysis, or wherefrom to step in synthesis;
nothing could be accurately understood,\\ nothing rigorously
demonstrated; and a full conviction never be obtained in
the pursuit of sciences.
}
\end{quotation}

The organisation of Darr\'e's text leaves something to be
desired, even apart from the quality of his English,
for which he frankly begs indulgence.  He does
not give explicit postulates, and his \lq\lq proof"
of the angle-sum theorem for triangles employs a
couple of hidden assumptions and a previous
result with a useless proof involving motion.
Nicholas Callan
later wrote what he describes as a 
revised and improved version of Darr\'e's text,
in which he assumes explicit Postulates
including a form of the parallel postulate\footnote{%
\lq\lq Postulate 4: A straight line that meets one of two
 parallels, may be produced until it meets the other."}.
His definition of parallel\footnote{\lq\lq having no divergency"}
is not terribly useful, and 
is possibly influenced by Legendre \cite{L},
two editions of which are in the College Library.
But it is clear that he broke with Darr\'e and Legendre
in making no attempt to do
without an axiom of parallels. Also, Legendre
tries to {\em prove} that all right angles
are equal, instead of just assuming it.  We suppose that
Callan \lq\lq went back to Euclid", to a large extent.

\section{Modern difficulties}

\subsection{Orientation}
Suppose that we draw a triangle and mark a small arrow-head on its 
boundary to indicate the sense in which we consider a moving point makes one 
complete circuit of the boundary. Visually it is very clear that there are 
two possibilities, one of which we name \textit{clockwise} and the other of 
which we name \textit{anticlockwise}. But how are we to put a definition or 
properties of that into our mathematical model? If you have drawn your 
triangle on a sheet of grease-proof paper you will see that the situation 
is reversed when looked at through the back of the paper, and that is one complication. We cannot put an arrow-head on each triangle boundary in the plane, 
so what can we do? Mathematicians have worked out a way to handle
this problem by placing an arrow-head on one boundary and improvising a 
method of transferring that to the boundary of every other triangle 
in the plane, so that, for example, we can say what is 
clockwise on every boundary. This concept is named {\em orientation} 
of a circuit on a triangle-boundary.

\subsection{Rotation}
A second awkward topic is the very familiar one of a {\em rotation} in 
the real world, or to put it more precisely,
\lq rotation about a given point, through a given angle'.  This involves the 
difficulty of orientation too. Mathematics has long had a formula for 
this in coordinate geometry, which uses 
trigonometry, but how can it be handled if it is introduced early on in 
pure geometry?

\subsection{How should these difficulties be handled?}

Our position is not that these are very difficult concepts but that they 
should wait until it is much easier to introduce them. The concepts 
should be made clear visually by diagrams.

For example orientation, which is rarely dealt with formally, can 
easily be handled (if someone wants to do that) by 
using the concept of sensed or signed area in 
coordinate geometry.

The syllabus documents as they stand do not include a formal treatment of 
rotations, but do mention them. If a formal treatment
 were to be added within the existing framework of five strands,
a rotation would be a type of function $P\rightarrow R(P)$. 
It would take points $P$ of the
plane as inputs, and give other points $Q=R(P)$ as outputs. 
The relationship
between the input point $P$, the centre of rotation $C$, and the output
point $Q$ would be described in terms of angles and congruence: 
the angles $POQ$ would all be congruent, and $|PC| = |QC|$.

\noindent
(A document in which a treatment of this kind is
presented as material for a group project may be downloaded
from either of the sites\\
{\Small \url{http://www.ucc.ie/en/euclid/edu_and_careers/projectmaths/} 
}
or\\
{\Small \url{http://archive.maths.nuim.ie/staff/aof/school.html}}. 
These sites also have a few other resources related
to school geometry, which have been submitted
for approval to the NCCA Project Maths coordinators.)

We note that the text-book of Hall and Stevens \cite{HS}, 
very commonly used in Ireland in the past, 
contained alternative proofs by rotation of some theorems, 
starting as early as Theorem 1. Teachers may be familiar with
these. It may seem attractive to use such proofs,
however, the point of training in deductive thinking is lost if proofs can pull in extra axioms out of the blue, and there are no
axioms about rotations in the present system.  Formal proofs
studied should remain within whatever logical framework
is laid down in the syllabus, and this in its turn
must be based upon 
a scholarly back-ground published treatment which provides a context for it.

\subsection{Angles and Rotations}
This is not to suggest that students should  not be made very 
familiar with the visual concepts of clockwise and 
anticlockwise rotations in the real world.
Of course they should.
According to the Primary Curriculum, students
are to \lq learn to recognise an angle
in terms of a rotation' \cite[p.75]{primary}.  This is
a bit ambiguous, but perhaps ok.  
In the formal material on geometry, angle is an undefined term,
an abstraction like point and line, 
so the question is: how is the student to be prepared 
for this, in the preliminary stage?  From this stage,
the student will bring some intuitive idea of what
an angle is.

It says in the syllabus document that to each angle
is associated a unique point called its vertex,
two rays starting at the vertex, called its arms,
and a piece of the plane called its inside.

It is not going to work very well if the student thinks
that an angle is \lq\lq a rotation".  This carries with
it some idea of motion, and this is not helpful
in studying the angles of, for example, a given triangle.
So it would be better, when talking informally
about rotations, to say that a rotation is something
that can be described in terms of an angle, rather than saying
that an angle is a rotation.  An angle is a specific
\lq\lq static" object, with vertex, arms and inside.

There are two things: the angle, and the number of degrees
in the angle (also known as the measure of the angle \cite{Barry}).
One is a geometrical object, the other is a real number.
These can be confused.  We suggest that it is a good
idea, in the preliminary informal work, to draw various
angles, point out the vertex, arms, inside of each,
and say that the number of degrees tells
us \lq\lq how big the opening is", \lq\lq how rapidly the arms
diverge", \lq\lq how much we have to {\bf rotate }one arm
about the vertex, in order to reach the other arm",
\lq\lq the amount of turning involved, if we first
face along one arm, and then turn and face along
the other",
and go onto discuss how much of a circle about
the vertex is inside the angle, the concepts of degree
(and radian, if desired) and the use of the protractor.  

\medskip
\noindent{\bf Note. }The references below
include some old books, long out of print
and probably not accessible to most teachers. Happily, most such 
out-of-copyright books may now be accessed and read 
online, thanks to the Google books initiative and other archiving
efforts such as archive.org.  We recommend that teachers, when time allows,
take advantage of these resources.

\end{document}